\documentclass[11pt]{article}
\usepackage[utf8]{inputenc}

\usepackage[margin=1in]{geometry}

\usepackage{comment}

\usepackage{subfig}

\usepackage{etex}
\usepackage{scrhack}
\usepackage{cmap}
\usepackage[english]{babel}
\usepackage[utf8]{inputenc}
\usepackage{caption}
\usepackage{epstopdf}
\usepackage{mathtools}
\usepackage{amsmath}
\usepackage{amssymb}
\usepackage{amsfonts}
\usepackage{graphicx}
\usepackage{tikz}
\usepackage[english,vlined,ruled]{algorithm2e}
\usepackage{float}
\usepackage{makeidx}
\usepackage{tabularx}
\usepackage{booktabs}
\usepackage{longtable}
\usepackage{tabu}
\usepackage{environ}
\usepackage{mdframed}
\usepackage{datetime}
\usepackage{enumerate}
\usepackage{hyperref}
\usepackage{cleveref}
\usepackage{todonotes}
\usepackage{appendix}
\usepackage{caption}
\usepackage{hyperref}

\DeclareMathOperator*{\argmax}{\arg\!\max}

\newcommand{\p}{\textsc{P}\xspace}
\newcommand{\np}{\textsc{NP}\xspace}

\newcommand{\forward}{(\( \Rightarrow \))\xspace}
\newcommand{\backward}{(\( \Leftarrow \))\xspace}

\newcommand{\define}[1]{\textit{#1}}

\newcommand{\osc}{\textsc{SC-Orientation}} %
\newcommand{\OSC}{SCO}
\newcommand{\tsat}{\textsc{3-SAT}}

\newcommand{\ladder}{\mathcal{L}}
\newcommand{\cltsat}{\textsc{Clause-linked Planar 3-SAT}}
\newcommand{\laddercycle}{\mathcal{LC}}

\newcommand{\GG}{\mathcal{G}}
\newcommand{\girth}{\operatorname{girth}}

\newcommand{\lemmaDHarePDH}{The sc-orientable distance hereditary graphs are exactly the strongly distance hereditary graphs.
	They can be recognized in polynomial time.}

\newcommand{\commentt}[1]{}

\newenvironment{prooff}{\begin{proof}}{\phantom{S} \end{proof}}

\usepackage{amsthm} %
\usepackage{xspace}
\usepackage{extarrows}
\usepackage[T1]{fontenc}
\usepackage{csquotes}

\usepackage{cite}

\usepackage{thm-restate}

\newtheorem{theorem}{Theorem}[section]
\newtheorem{lemma}[theorem]{Lemma}

\newcommand{\X}[1]{}

\title{Make a graph singly connected by edge orientations}

\author{
Tim A.~Hartmann\footnote{CISPA Helmholtz Center for Information Security, Saarbrücken, Germany \texttt{tim.hartmann@cispa.de}}
\and Komal Muluk\footnote{RWTH Aachen University, Germany. \texttt{muluk@algo.rwth-aachen.de}}
}

\date{\today}

\begin{document}

\maketitle

\begin{abstract}
\begin{sloppypar}
	A \emph{directed} graph $D$ is \emph{singly connected} if for every ordered  pair of vertices $(s,t)$, there is at most one path from $s$ to $t$ in $D$. 
	Graph orientation problems ask, given an \emph{undirected} graph $G$,
	to find an orientation of the edges such that the resultant \emph{directed} graph $D$ has a certain property.
	In this work, we study the graph orientation problem where the desired property is that $D$ is singly connected.
	Our main result concerns graphs of a fixed girth $g$ and coloring number $c$.
	For every $g,c\geq 3$, the problem restricted to instances of girth $g$ and coloring number $c$, is either~\np-complete or in~\p.
	As further algorithmic results, we show that the problem is \np-hard on planar graphs and polynomial time solvable distance-hereditary graphs.
\end{sloppypar} 
\end{abstract}

\section{Introduction}
\label{section:introduction}
Graph orientation problems are well-known class of problems in which given an undirected graph, the task is to decide whether we can assign
a direction to each edge such that a particular desired property is satisfied on the resulting directed graph.
One can ask questions like, for an undirected graph, can we assign directions to all edges such that the digraph thus formed is acyclic, or it is strongly connected, or it contains a directed Euler circuit, and many more such interesting properties. These problems are important because of the emerging field of networks designing; networks are modeled by directed graphs. %
Consider the problem of converting an undirected graph to a directed acyclic graph by means of such an orientation of edges. All simple graphs can be converted into directed acyclic graphs by simply considering a DFS tree on the graph and directing all edges from the ancestor to the descendant. This trivially solves the decision version of the graph orientation problem for all simple graphs.

In this work we consider `singly-connected' as the desired property of the graph orientation property.
A directed graph $D$ is \emph{singly connected} if for every ordered pair of vertices $(s,t)$, there is at most one path from $s$ to $t$ in $D$.  %
For an undirected graph $G$, an $E(G)$-orientation is a mapping $\sigma$ that maps each edge $\{u,v\} \in E(G)$ to either of its endpoints $u$ or $v$, which results into a directed graph $G_\sigma$.
An $E(G)$-orientation $\sigma$ is an \define{sc-orientation} if $G_{\sigma}$ is singly connected. 
A graph $G$ is said to be \emph{sc-orientable} if there exists an sc-orientation of $G$.
Hence we study the following problem.
\\
\fbox
{\begin{minipage}{42em}
		\textbf{Problem:} \osc\ (\OSC) \\
		\textbf{Input:} A simple undirected graph $G=(V,E)$.   \\
		\textbf{Question:} Is $G$ sc-orientable?
\end{minipage}}\\

There has been several studies around the singly connected graphs in the literature. The problem of testing whether a graph is singly connected or not was introduced in an exercise in Cormen et al.\ in \cite{Cormen}. Buchsbaum and Carlisle, in 1993, gave an algorithm running in time $O(n^2)$ \cite{Buchsbaum}. There has been other attempts to improve this result, for example \cite{Khuller,DBLP:journals/ipl/DietzfelbingerJ15}.
For our \osc\ problem, this means that the edge orientation $\sigma$
that makes $G_\sigma$ singly connected serves as an \np-certificate.
Thus \OSC\ is in \np. 

Das et al.~\cite{DasKMMPS2022} studied the vertex deletion version of the problem. They asked, given a \emph{directed} graph $D$, does there exists a set of $k$ vertices whose deletion renders a singly connected graph.
Alongside, they pointed a close link of the \textsc{Singly Connected Vertex Deletion (SCVD)} problem to the classical \textsc{Feedback Vertex Set (FVS)} problem. They argued how SCVD is \emph{a} directed counterpart of FVS, which makes it theoretically important. Our problem is a natural follow-up question in this direction which looks at the graph orientation version.

In the next section we prove that if a graph has an sc-orientation, then it also has an acyclic sc-orientation.
So, in addition for a graph to orient to an acyclic digraph, \osc\ also demands the orientation to a singly connected digraph. As discussed earlier, it is trivial to find an acyclic orientation of the graph. However, in contrast to that, we prove, the additional singly connected condition makes it harder to even decide whether such an orientation exits or not.

\subsection{Our Results}

First, we prove that \osc\ is \np-complete even for planar graphs. 
Thus the hardness also holds for all graphs restricted to girth at most 4.
Complementing this result, we show that planar graphs of girth at least 5 are always sc-orientable, hence, the decision problem is trivial.

We aim to show analogous dichotomy for general graphs.
As an upper bound, we show, for any graph with girth at least twice the chromatic number, the graph is sc-orientable.
Hence the problem is trivially solvable for such inputs.
To account for this, we study \osc\ restricted to the graph class $\GG_{g,c}$, the class of graphs that have girth $g$ and chromatic number $c$, for some positive integers $g,c$.
Our main result is that, for each pair $g,c\geq 3$, \osc\ restricted to $\GG_{g,c}$ %
is either \np-complete or trivially solvable.
We do not pinpoint the exact boundary for the values of $g$ and $c$ where the transformation occurs.
The approach is to compile an \np-hardness proof that solely relies upon a single no-instance of the given girth $g$ and chromatic number~$c$.

As further results, we provide polynomial time algorithms of \OSC\ for a special graph class such as distance-hereditary graphs.
Here, we give a concise definition of the yes-instances by forbidden subgraphs.

From the technical sides, we provide two simplification of the problem.
First, one can restrict the search for an sc-orientation to the orientations that also avoid directed cycles.
Second, we show how we can assume that the input is triangle free.

\medskip

In \Cref{section:preliminaries}, we begin with the Preliminaries and the structural results.
\Cref{section:planar} shows \np-completeness of \osc\ on planar graphs.
\Cref{section:girth} contains our dichotomy result.
Then, \Cref{section:perfect} considers perfect and distance-hereditary graphs.
We conclude with \Cref{section:conclusion}.

\section{Notation and Technical Preliminaries}
\label{section:preliminaries}

	\def\boxSize{2.5cm}
	\def\sideSep{0.2cm}
	\def\sideDist{0.25cm}
	
\begin{figure}[t]
\begin{center}
\captionsetup{position=b}
\begin{tikzpicture}[scale=1, node/.style={draw, circle, fill, minimum size=.15cm,inner sep=0pt}]
\begin{scope}[local bounding box=boxA]

\node[node] (a1) at (0,0) {};
\node[node] (a2) at (1,0) {};
\node[node] (b1) at (0,1) {};
\node[node] (b2) at (1,1) {};
\foreach \from/\to/\lbl in {a1/a2,b1/b2,a1/b2,a1/b1,a2/b2}
    \draw (\from) -- (\to);

\end{scope}
\begin{scope}[local bounding box=boxA,xshift=2cm]

\node[node] (a1) at (0,0) {};
\node[node] (a2) at (1,0) {};
\node[node] (b1) at (0,1) {};
\node[node] (b2) at (1,1) {};
\node[node] (c) at (0.5,1.8) {};
\foreach \from/\to/\lbl in {a1/a2,b1/b2,a1/b1,a2/b2,b1/c,b2/c}
    \draw (\from) -- (\to);

\end{scope}
\begin{scope}[local bounding box=boxA,xshift=4cm]

\node[node, label=below:$a_1$] (a1) at (0,0) {};
\node[node, label=below:$b_1$] (a2) at (1,0) {};
\node[node] (b1) at (0,1) {};
\node[node] (b2) at (1,1) {};
\node[node, label=above:$a_2$] (c1) at (0,2) {};
\node[node, label=above:$b_2$] (c2) at (1,2) {};
\foreach \from/\to/\lbl in {a1/a2,b1/b2,a1/b1,a2/b2,c1/c2,b1/c1,b2/c2}
    \draw (\from) -- (\to);
\end{scope}
\begin{scope}[local bounding box=boxA,xshift=7.5cm]

\node[node] (a1) at (0,0) {};
\node[node] (b0) at (-1.5,1) {};
\node[node] (b1) at (-0.5,1) {};
\node[node] (b2) at (0.5,1) {};
\node[node] (b3) at (1.5,1) {};
\foreach \from/\to/\lbl in {a1/b0,a1/b1,a1/b2,a1/b3,b0/b1,b1/b2,b2/b3}
    \draw (\from) -- (\to);

\end{scope}
\begin{scope}[local bounding box=boxA,xshift=10cm]

\node[node] (a1) at (0,0) {};
\node[node] (a2) at (0.7,1) {};
\node[draw=none,below =0.5cm of a2] (l1) {$H_1$};
\node[node] (a3) at (1.4,0) {};
\node[draw=none,above =0.5cm of a3] (l1) {$H_2$};
\node[node] (a4) at (0.7*3,1) {};
\node[draw=none,below =0.5cm of a4] (l1) {$H_3$};
\node[node] (a5) at (0.7*4,0) {};
\node[draw=none,above =0.5cm of a5] (l1) {$H_4$};
\node[node] (a6) at (0.7*5,1) {};

\foreach \from/\to/\lbl in {a1/a2,a2/a3,a3/a4,a4/a5,a5/a6}
    \draw (\from) -- (\to);
\draw[fill=black,opacity=0.2](0,0) -- (0.7,1) -- (1.4,2) -- (0.7*3,1) --(0.7*4,2) --(0.7*5,1) --(0.7*4,0) --(0.7*3,-1) --(0.7*2,0)-- (0.7,-1) -- (0,0); 
\end{scope}
\end{tikzpicture}
\end{center}
\caption{
From left to right: a diamond, a house, a domino (a $G_{2,3}$), a gem, and a gadget for the construction of \Cref{{lemma:no:iff:coupling}}.
Non-sc-orientable are the diamond, the house and hence also the gem.
The domino is sc-orientable.}
\label{figure:small:graphs}
\end{figure}
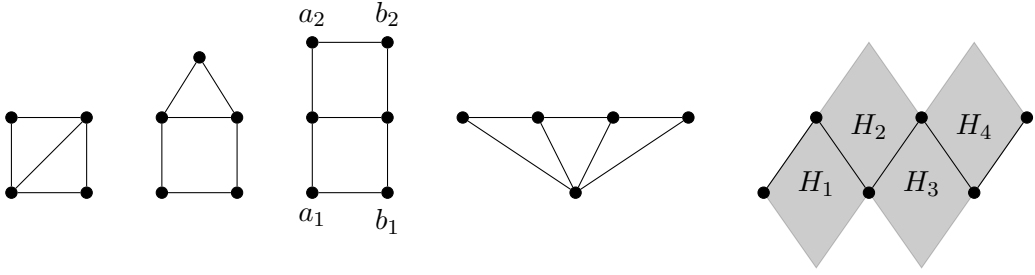

For a graph $G$, let $V(G)$ be its vertex set,
and $E(G)$ be its set of edges. Elements in $V(G)$ are called  nodes or vertices.
A $k$-coloring of a graph $G$, for $k \in \mathbb{N}$, is a mapping $f \colon V(G) \to \{1, \dots, k\}$ such that $f(u)\neq f(v)$ if $\{u,v\}\in E(G)$.
The coloring number, $\chi(G)$, of graph $G$ is the smallest integer $k$ such that $G$ has a $k$-coloring.
A graph has girth $g$ if the length of the shortest cycle in the graph is $g$. By definition, a graph without any cycles has girth $\infty$.
Given a path $P=v_1, v_2, \dots, v_n$, we denote the length of a path to be the total number of vertices in the path.
Throughout the paper we consider different small graphs, typically as forbidden induced subgraphs, drawn in \cref{figure:small:graphs}.
We denote the grid graph of width $x$ and height $y$ as $G_{x,y}$.
For further reading related to the concepts in graph theory we refer reader to the book by Diestel~\cite{DBLP:books/daglib/0030488}.

\medskip

The following technical observations are useful for our later results. The upcoming lemma shows when looking for an sc-orientation of a graph, it suffices to look for an sc-orientation that additionally avoids cycles of length greater than $3$.

\newcommand{\lemmaCycle}{
	A graph $G$ is sc-orientable
	if and only if
	there is an sc-orientation $\sigma$ of $G$ where $G_\sigma$ contains no directed cycle of length $\geq 4$.
}
\begin{lemma}%
	\label{lemma:no:directed:cycles}
	\lemmaCycle
\end{lemma}
\begin{prooff}
	The reverse direction is clear from the definition of sc-orientable graphs.
	
	For the forward direction, since $G$ is sc-orientable, let $\sigma$ be any sc-orientation of $G$. If $G_\sigma$ does not contain any directed cycle of size $\geq 4$, then we are done. Otherwise let $v_1,v_2,\dots,v_k, v_1$,  $k \geq 4$  be a directed cycle $C$ in $G_\sigma$. Let us denote the set of vertices of $C$ by $V(C)$. Now let the orientation $\sigma'$ be same as $\sigma$ except the orientation of the edges $(v_1, v_2)$ and $(v_3, v_4)$ reversed, i.e., unlike $G_\sigma$, graph $G_{\sigma'}$ contains edges $(v_2, v_1)$ and $(v_4,v_3)$.
	Let us refer to the set of newly oriented edges $\{(v_2, v_1), (v_4,v_3)\}$ by $E'$. 
	We claim that $\sigma'$ is an sc-orientation and it does not create any new directed cycle in $G_{\sigma'}$ apart from those present in $G_{\sigma}$.
	
	Let us first prove that $G_{\sigma'}$ is singly connected.
	Assume not, and $s, t\in V(G_{\sigma'})$ be such that there are two vertex disjoint paths from $s$ to $t$, $P_1'$ and $P_2'$, in $G_{\sigma'}$.
	At least one of $P_1',P_2'$ must use an edge from $E'$, otherwise $P_1'$ and $P_2'$ are present in $G_\sigma$, contradicting to the fact that $G_\sigma$ is an sc-orientation.
	We consider the following cases:
	\begin{itemize}
		\item \textit{Case 1:}
		Exactly one of $P_1', P_2'$, say $P_1'$, contains vertices from $V(C)$.
		And $P_2'$ does not intersect with the vertex set $V(C)$. This implies there exists a path $P_2=P_2'$ in $G_{\sigma}$. Now while tracing the path $P_1'$ from $s$ to $t$, consider the first vertex and the last vertex on $P_1'$ that belongs to $V(C)$. Let the vertices be $v_i$ and $v_j$ respectively; there exists a path $P$ joining $v_i$ to $v_j$ in the cycle $C$ of $G_\sigma$. Then by replacing the path between $v_i$ to $v_j$ in $P_1'$ by $P$, we can obtain a path $P_1$ from $s$ to $t$ in $G_\sigma$. Note that $P_1\neq P_2$, as $P_1$ contains vertices from $V(C)$. Thus we get two disjoint paths in $G_\sigma$, contradicting to the fact that $G_\sigma$ is singly connected.
		
		\item \textit{Case 2:}
		Both $P_1',P_2'$ intersect $V(C)$.
		Assume that either $s\notin V(C)$ or $t\notin V(C)$. Then without loss of generality, consider $s \notin V(C)$. 
		Tracing $P_1'$ from $s$ to $t$, $v_i$ be the first vertex that belong to $V(C)$. 
		Similarly, while tracing $P_2'$ from $s$ to $t$, $v_j$ be the first vertex that belong to $V(C)$. 
		If $v_i=v_j$, we get two paths from $s$ to $v_i$ in $G_\sigma$, a contradiction to $\sigma$ being an sc-orientation. 
		Thus $v_i\neq v_j$; there exists a path $P$ from $v_i$ to $v_j$ in the cycle $C$. 
		In this scenario, we get two paths from $v_i$ to $v_j$ in $G_\sigma$, one the subpath $P_2$ of $P_2'$ from $s$ to $v_j$, and another by concatenating the subpath $s, \dots, v_i$ of $P_1'$ with $P$. 
		Both are disjoint paths from $s$ to $v_j$ in $G_\sigma$, giving us a contradiction to $\sigma$ being an sc-orientation. 
		Similarly, if we assume that $t\notin V(C)$, we attain the same contradiction. 
		
		\item \textit{Case 3:}
		Both $s,t\in V(C)$. Note that at least one of $P_1', P_2'$, say $P_1'$, contains at least one edge $(u,v)$ such that $u\in V(C)$ and $v\notin V(C)$.  
		Then while tracing $P_1'$ from $s$ to $t$, consider the last vertex $v_i$ such that $P_1'$ contains an edge $(v_i, w)$, where $ v_i\in V(C) \text{ and } w\notin V(C)$. 
		The vertex $v_j\in V(C)$ be such that $P_1'$ contains a subpath $v_i, w_1, \dots, w_n, v_j$ where $w_i\notin V(C)$ for $i\in [n]$. 
		Then we obtain two paths from $v_i$ to $v_j$ in $G_\sigma$, namely $P_1=v_i, w_1, \dots, w_n, v_j$, and $P_2$ can be traversed on the cycle $C$. 
		Thus giving two disjoints paths from $v_i$ to $v_j$ in $G_\sigma$, a contradiction.  
	\end{itemize}
	
	So far we have proved that $\sigma'$ is an sc-orientation. To prove our claim we only need to show that the newly oriented edges from $E'$ do not create new directed cycles in the graph $G_{\sigma'}$.
	Assume the contrary. If there is a new cycle $C'$ in $G_{\sigma'}$ that was not present in $G_\sigma$, then this cycle must include an edge from $E'$. If $C'$ contains only one edge from $E'$, say $(v_2, v_1)$, then there is a path from $v_1$ to $v_2$ in $G_\sigma$ which does not contain the edge $(v_1,v_2)$, thus creating two disjoint paths in $G_\sigma$, a contradiction. 
	
	\begin{sloppypar}
	Now if $C'$ contains both the edges from $E'$, namely $(v_2,v_1)$ and $(v_4, v_3)$, then let $C'=v_2, v_1, w_1, w_2, \dots, w_i, v_4, v_3, w_{i+1}, \dots, w_j, v_2$.
	Note that the vertices $w_1, w_2, \dots, w_j\notin \{v_1, v_2, v_3, v_4\}$. Thus in this case, we get two paths from $v_1$ to $v_4$, namely $v_1, v_2, v_3, v_4$ and $v_1, w_1, w_2, \dots, w_i, v_4$ in $G_\sigma$. 
	This concludes the proof because each time we encounter a directed cycle of size $\geq 4$, we can reduce the number of directed cycles by at least 1.
	\end{sloppypar}
\end{prooff}

Further, we can remove triangles from the input graph unless it contains a diamond or a house as an induced subgraph. It is easy to see, if a graph $G$ contains a diamond or a house graph as an induced subgraph, then $G$ is not sc-orientable.

\newcommand{\lemmaTriangle}{
	Let a (diamond, house)-free graph $G$ contain a triangle $uvw$.
	Then $G$ has an sc-orientation, if and only if $G'$, resulting from $G$ by contracting the edges in $uvw$, has an sc-orientation.
}
\begin{lemma}%
	\label{lemma:triangle:shrink}
	\lemmaTriangle
\end{lemma}
\begin{prooff}
	We contract the triangle $uvw$ in $G$ to a single vertex in $G'$. Let $z$ be that single vertex in $G'$.
	
	($\Rightarrow$)
	Let $\sigma$ be an sc-orientation of $G$; the triangle $uvw$ obtains a cyclic orientation in $G_\sigma$.
	We show that then $\sigma'$ is an sc-orientation of $G'$ where $\sigma'$ is defined as follows:
	\[
	\sigma'(\{x,y\}) =
	\begin{cases}
		\sigma(\{x,y\}), & \text{if $x,y\neq z$}, \\
		\sigma(\{x,z'\}), & \text{if $y=z$ and $\{x,z'\} \in E(G)$ where $z'\in \{u,v,w\}$}. \\
	\end{cases}
	\]
	
	Note that since $G$ is a diamond-free graph, an edge $xz\in E(G')$ implies there is a unique $y\in \{u,v,w\}$ such that $xy\in E(G)$. Now, for the sake of contradiction, assume that $\sigma'$ is not an sc-orientation, and there are two internally vertex disjoint $s,t$-paths $P'_1 = x_0, x_1, x_2,\dots, x_{m-1}, x_m$, and $P'_2= y_0,y_1, y_2, \dots, y_{n-1}, y_n$
	for some $s=x_0=y_0$, $t=x_m=y_n \in V(G')$.
	Then either $z\in \{s,t\}$, or $z$ belongs to exactly one of $P_1', P_2'$; if $z$ does not belong to either of the paths, then $P_1', P_2'$ are present in $G_\sigma$ giving us a contradiction. If $z\in \{s,t\}$ (say $z=s$),
	then there exist arcs $(a,x_1), (b,y_1)$ in $G_\sigma$ where $a,b\in \{u,v,w\}$.
	We thus obtain two a,t-paths, or two b,t-paths depending on the orientation of the edges of triangle $uvw$ in $G_\sigma$.
	
	Now, if $z$ belongs to only one path, say $z$ belongs to $P_1'$. Assume $z=x_i$. Then there are vertices $a, b\in \{u,v,w\}$, such that the arcs $(x_{i-1}, a),(b,x_{i+1})$ belongs to $G_\sigma$. Let $P$ be an a,b-path obtained by traversing the edges of triangle $uvw$ of $G_\sigma$. Then $P_1=s, \dots, x_{i-1}, P, x_{i+1}, \dots,t$, and $P_2\coloneqq P_2'$ are disjoint $s,t$-paths in $G_\sigma$, hence a contradiction.
	
	($\Leftarrow$)
	Let $\sigma'$ be an sc-orientation of $G'$.
	With the help of \Cref{lemma:no:directed:cycles} we assume that $G_{\sigma'}'$ does not contain a directed cycle of length at least $4$.
	We then show that $\sigma$ is an sc-orientation of $G$ where $\sigma$ is obtained from $\sigma'$, and for the edges of the triangle $uvw$ it orients the edges to obtain a cyclic triangle. Assume the contrary, there are two internally vertex disjoint $s,t$-paths $P_1,P_2$ in $G_\sigma$ for some $s,t \in V(G)$. At least one of $P_1, P_2$ contains vertices from the set $\{u,v,w\}$. Then we distinguish among the following cases:
	\begin{itemize}
		\item
		Consider that at most one path of $P_1,P_2$ goes through the vertices of the  triangle $uvw$. Say only $P_1 = s, v_1,\dots, v_k, t$ does so, and let the first and last vertices in $P_1 \cap \{u,v,w\}$ be $v_i$ and $v_j$ respectively.
		Note that $s,t \notin \{u,v,w\}$. 
		Then let $P_2' \coloneqq P_2$ and $P_1'$ be same as $P_1$ but where the subpath $v_i,\dots,v_j$ is replaced by the vertex $z$ in $G_{\sigma'}'$.
		Then $P_1',P_2'$ are disjoint $s,t$-paths in $G_{\sigma'}'$, a contradiction.
		\item
		Now consider if both paths $P_1, P_2$ traverse through the vertices of the triangle $uvw$. However, we assume that $s,t\notin \{u,v,w\}$. Then $v_{i_1}, v_{i_2}$ be the corresponding first vertices of $P_1, P_2$ that belongs to $\{u,v,w\}$. Path $P_1'$ be the subpath of $P_1$ from $s$ to $v_{i_1}(=z$ in $G')$, and $P_2'$ be the subpath of $P_2$ from $s$ to $v_{i_2}(=z$ in $G')$. We obtain two disjoint $s,z$-paths in $G_{\sigma'}'$.
		\item
		Given that $P_1, P_2$ both traverse through the triangle, we further assume that exactly one of $s, t$(say $s$) belongs to $\{u,v,w\}$. Here $t\notin \{u,v,w\}$. 
		Let $v_{j_1}$ be the last vertex of $P_1$ such that $v_{j_1}\in \{u,v,w\}$. Similarly, $v_{j_2}$ be the last vertex of $P_2$ such that $v_{j_2}\in \{u,v,w\}$. Then $P_1'= v_{j_1}(=z \text{ in }G'), \dots, t$, and $P_2'= v_{j_2}(=z \text{ in } G'), \dots, t$ are disjoint $z,t$-paths in $G_{\sigma'}'$.
		Thus we get a contradiction.
		\item
		It remains to consider that $\{s,t\} \subseteq \{u,v,w\}$.
		Since $uvw$ form a directed cycle in $G_\sigma$,
		at least one of the paths $P_1,P_2$, say $P_1$, visits vertices not in $\{u,v,w\}$.
		Let $v_1,v_2,\dots,v_i$ be a subpath of $P_1$ where $v_1,v_i \in \{u,v,w\}$ and $v_2,\dots,v_{i-1} \notin \{u,v,w\}$.
		Because $G$ is diamond- and house-free, $i\geq 5$.
		Then  $G_{\sigma'}'$ contains a directed cycle $z,v_2,\dots,v_{i-1},z$ of length at least $4$. A contradiction to our initial assumption that $G_{\sigma'}'$ does not contain any directed cycle of length at least 4. 
	\end{itemize}
	Since all cases led to a contradiction, $G_{\sigma}$ must be singly connected.
\end{prooff}

\section{\np-hardness on Planar Graphs}
\label{section:planar}

This section derives \np-hardness of \osc\ on planar graphs.
As observed in the \cref{section:introduction}, \osc\ is in \np.
To show hardness, we give a reduction from a variant (to be shortly defined) of the \textsc{Planar \tsat} problem. Given a \tsat\ formula $\phi$, a \emph{variable-clause graph} is a graph obtained by taking variables and clauses representatives as the vertex set, and the edge set contains an edge between a variable-clause pair if the variable appears in the clause. We use the variant of \tsat\  where there is a planar embedding of the variable-clause graph, in which there are some additional edges (only connecting clause nodes) which form a cycle traversing through all clause nodes.  
We denote this cycle as \emph{clause-cycle}, and the problem as \cltsat.
Kratochv{\'{\i}}l et al.\ proved the \np-hardness of \cltsat~\cite{DBLP:journals/siamdm/KratochvilLN91}.\\

\fbox
{\begin{minipage}{42em}
		\textbf{Problem: }\cltsat\\
		\textbf{Input:} A boolean formula $\phi$ with clauses $C = \{c_1,\dots,c_m\}$ and variables $X=\{v_1, \dots, v_n\}$,
		and a planar embedding of the variable-clause graph
		together with a clause-cycle.
		\\
		\textbf{Question:} Is $\phi$ satisfiable?
\end{minipage}}\\

\begin{theorem}%
	\label{lemma:np:planar}
	\osc\ is \np-complete even for planar graphs.
\end{theorem}
\begin{prooff}
	Our idea is to take the given planar embedding, %
	and use this as a backbone for our construction.
	Let $\eta$ be an ordering of clause nodes
	such that the clause-cycle traverses through $C$ as per the ordering $\eta$.
	Let the given planar embedding be $\Pi$, and the underlying graph be $G_{\eta}(\phi)$. 
	For the reduction, %
	we define a clause gadget and a variable gadget to replace the clause nodes and variable nodes in $\Pi$. 
	To begin with, let us start with the following interesting graph structure:
	
	\medskip
	\noindent\textbf{Ladder.}
	The domino (the $2\times 3$-grid graph, see \cref{figure:small:graphs}),
	has exactly two possible sc-orientations.
	If we fix the orientation of one edge, all orientations of the domino are fixed up to singly-connectedness
	(every vertex becomes either a sink or a source for the domino).
	Specifically, we are interested in the following properties
	regarding the bottom edge $\{a_1,b_1\}$ and the top edge $\{a_2,b_2\}$:
	\begin{enumerate}[{(1)}]
		\item \sloppy
		$(a_1,b_1), (a_2,b_2)$ are \define{coupled},
		that is, for every sc-orientation $\sigma$,
		either $\sigma(\{a_x,b_x\}) = b_x$ for $x\in\{1,2\}$,
		or $\sigma(\{a_x,b_x\}) = a_x$  for $x\in\{1,2\}$,\label{prop1}
		\item
		there is an sc-orientation without any directed path from $\{a_1,b_1\}$ to $\{a_2,b_2\}$.
	\end{enumerate}
	
	Since the directions of coupling edges matter, we refer to them as ordered pairs, for example, edge $\{a,b\}\in V(G)$ becomes $(a,b)$ or $(b,a)$.
	The coupling property is transitive,
	in the sense that, when we take 2 dominos and identify the top edge of one to the bottom edge of another (forming a $2 \times 5$-grid), the very top and the very bottom edges are also coupled. Further extrapolation can be made to form a \emph{ladder} (a $2\times n$-grid).
	Additionally, the bottom and the middle edges of a domino are in the reverse direction in every sc-orientation. Hence they are the \emph{reverse-coupled} edges.
	The coupling property of dominos helps us attach as many $4$-cycles as we want to a ladder without affecting the two possible sc-orientations of the structure. Let us call the ladder with the additional $4$-cycles an \emph{extended ladder} $(\ladder)$. For our reference, let us call the bottom of the ladder the $0^{th}$-step.
	Let us call the set of edges which are coupled with the $0^{th}$-step by \emph{even edges}, and the set of reverse-coupled edges by \emph{odd edges}.
	All in all, if we fix the orientation of the $0^{th}$-step, the orientation of the whole ladder is fixed up to an sc-orientation. Let us construct a large enough extended ladder $\ladder_{U}$, we call this the universal ladder.
	\medskip

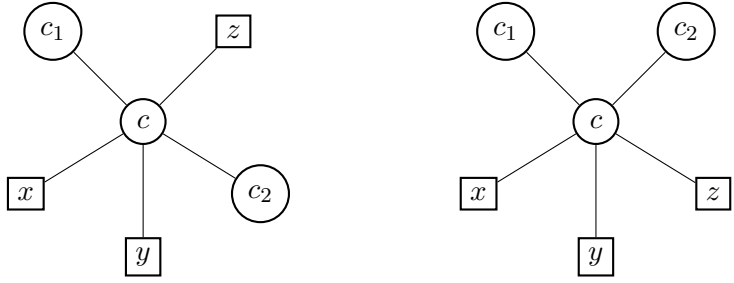
\begin{figure}[t]
    \begin{center}
    \def\prop{thick}
        \subfloat{%
        \begin{tikzpicture}[scale=1.2, >=latex]
            \node at (0,0) [circle,draw, \prop] (c) {$c$};
            \node at (-1,1) [circle,draw,\prop] (c_1) {$c_1$};
            \node at (1,1) [rectangle,draw,\prop] (z) {$z$};
            \node at (-1.3,-0.8) [rectangle,draw,\prop] (x) {$x$};
            \node at (0,-1.5) [rectangle,draw,\prop] (y) {$y$};
            \node at (1.3,-0.8) [circle,draw,\prop] (c_2) {$c_2$};
            \draw (c) -- (c_1);
            \draw (c) -- (c_2);
            \draw (c) -- (x);
            \draw (c) -- (y);
            \draw (c) -- (z);
        \end{tikzpicture}
        }
\hspace{20mm}
        \subfloat{%
        \begin{tikzpicture}[scale=1.2, >=latex]
            \node at (0,0) [circle,draw,\prop] (c) {$c$};
            \node at (-1,1) [circle,draw,\prop] (c_1) {$c_1$};
            \node at (1,1) [circle,draw,\prop] (c_2) {$c_2$};
            \node at (-1.3,-0.8) [rectangle,draw,\prop] (x) {$x$};
            \node at (0,-1.5) [rectangle,draw,\prop] (y) {$y$};
            \node at (1.3,-0.8) [rectangle,draw,\prop] (z) {$z$};
            \draw (c) -- (c_1);
            \draw (c) -- (c_2);
            \draw (c) -- (x);
            \draw (c) -- (y);
            \draw (c) -- (z);
        \end{tikzpicture}
        }
        
    \end{center}
    \caption{Type 1 (left) and type 2 (right) of clause nodes according to their neighbourhoods in the planar embedding $\Pi$ of $G_\eta(\phi)$.}
    \label{types_of_clause_nodes}
\end{figure}

\begin{figure}[t]
    \begin{center}
        \subfloat{%
        \begin{tikzpicture}[scale=0.8,>=latex]
    \def\vsize{0.6mm}
    \def\col{black!20}
    \draw [line width=2mm, \col] (1.5,3.5) -- (3.5,1.5);    
    \draw [line width=2mm, \col] (3.5,-1.5) -- (1.5,-3.5);    
    \draw [line width=2mm, \col] (-1.5,-3.5) -- (-3.5,-1.5);   
    \draw [line width=1mm] (-3.5,1.5) -- (-3.5,-1.5);
    \draw (-3.5,1.5) node (o1) {};\filldraw (o1) circle (\vsize);
    \draw (-1.5,3.5) node (o2) {};\filldraw (o2) circle (\vsize);
    \draw (3.5,1.5) node (o3) {};\filldraw (o3) circle (\vsize);
    \draw (1.5,3.5) node (o4) {};\filldraw (o4) circle (\vsize);
    \draw (3.5,-1.5) node (o5) {};\filldraw (o5) circle (\vsize);
    \draw (1.5,-3.5) node (o6) {};\filldraw (o6) circle (\vsize);
    \draw (-1.5,-3.5) node (o7) {};\filldraw (o7) circle (\vsize);
    \draw (-3.5,-1.5) node (o8) {};\filldraw (o8) circle (\vsize);
    \draw[black,thick]   (-3.5,1.5) -- (-1.5,3.5)  -- (1.5,3.5)  -- (3.5,1.5) -- (3.5,-1.5) -- (1.5,-3.5)  -- (-1.5,-3.5)  -- (-3.5,-1.5) -- (-3.5,1.5);
    \node[scale=1] at (-3.8,0.2) {\color{black}$e_1^c$};
    \node[scale=1] at (-2.7,-2.7) {\color{black}$e_x^c$};
    \node[scale=1] at (0,-3.8) {\color{black}$e_2^c$};
    \node[scale=1] at (2.7,-2.7) {\color{black}$e_y^c$};
    \node[scale=1] at (3.8,0.2) {\color{black}$e_3^c$};
    \node[scale=1] at (2.7,2.7) {\color{black}$e_z^c$};
    \node[scale=1] at (0,3.8) {\color{black}$e_4^c$};
    \node[scale=1] at (-2.7,2.7) {\color{black}$e_5^c$};
    \node[scale=1] at (-3.95,-0.3) {\color{black}$in_c$};
    \node[scale=1] at (4.05,-0.3) {\color{black}$out_c$};
    \draw (-5,1.5) node (in_c_11) {};\filldraw (in_c_11) circle (\vsize);
    \draw (-5,-1.5) node (in_c_12) {};\filldraw (in_c_12) circle (\vsize);
    \draw (-3,1.5) node (n1) {}; \filldraw (n1) circle (\vsize);
    \draw (-2.5,1.5) node (n2) {};\filldraw (n2) circle (\vsize);
    \draw (-2,1.5) node (n3) {};\filldraw (n3) circle (\vsize);
    \draw (-1.5,1.5) node (n4) {};\filldraw (n4) circle (\vsize);
    \draw (1.5,1.5) node (n5) {};\filldraw (n5) circle (\vsize);
    \draw (2,1.5) node (n6) {};\filldraw (n6) circle (\vsize);
    \draw (2.5,1.5) node (n7) {};\filldraw (n7) circle (\vsize);
    \draw (3,1.5) node (n8) {};\filldraw (n8) circle (\vsize);
    \draw (-3,-1.5) node (n9) {};\filldraw (n9) circle (\vsize);
    \draw (-2.5,-1.5) node (n10) {};\filldraw (n10) circle (\vsize);
    \draw (-2,-1.5) node (n11) {};\filldraw (n11) circle (\vsize);
    \draw (-1.5,-1.5) node (n12) {};\filldraw (n12) circle (\vsize);
    \draw (1.5,-1.5) node (n13) {};\filldraw (n13) circle (\vsize);
    \draw (2,-1.5) node (n14) {};\filldraw (n14) circle (\vsize);
    \draw (2.5,-1.5) node (n15) {};\filldraw (n15) circle (\vsize);
    \draw (3,-1.5) node (n16) {};\filldraw (n16) circle (\vsize);
    \draw (5,1.5) node (out_c_21) {};\filldraw (out_c_21) circle (\vsize);
    \draw (5,-1.5) node (out_c_22) {};\filldraw (out_c_22) circle (\vsize);
    \def\colorladder{black}
    \draw[black,thick] (-5,1.5) -- (-3.5,1.5) -- (-3,1.5) -- (-2.5,1.5) -- (-2,1.5) -- (-1.5,1.5) -- (1.5,1.5) -- (2,1.5) --(2.5,1.5) --(3,1.5) -- (3.5,1.5) -- (5,1.5);
    \draw[black,thick] (-5,-1.5) -- (-3.5,-1.5) -- (-3,-1.5) -- (-2.5,-1.5) -- (-2,-1.5) -- (-1.5,-1.5) -- (1.5,-1.5) -- (2,-1.5) --(2.5,-1.5) --(3,-1.5) -- (3.5,-1.5) -- (5,-1.5);
    \draw[draw, thick] (-1.5, 3.5)-- (-1.5,1.5) -- (-1.5,-1.5) -- (-1.5,-3.5);
    \draw[draw, thick] (1.5, 3.5)-- (1.5,1.5) -- (1.5,-1.5) -- (1.5,-3.5);
    \draw[thick] (-3,1.5) --(-3,-1.5); 
    \draw[thick] (-2.5,1.5) --(-2.5,-1.5); 
    \draw[thick] (-2,1.5) --(-2,-1.5); 
    \draw[thick] (3,1.5) --(3,-1.5); 
    \draw[thick] (2.5,1.5) --(2.5,-1.5); 
    \draw[thick] (2,1.5) --(2,-1.5);
    \draw[thick] (-3.5,1.5) ..controls(-2.7,2) .. (-2,1.5);
    \end{tikzpicture}
        }
        \subfloat{%
        \begin{tikzpicture}[scale=0.8]
    \def\vsize{0.6mm}
    \def\col{black!20}
    \draw [line width=2mm, \col] (1.5,3.5) -- (3.5,1.5);    
    \draw [line width=2mm, \col] (3.5,-1.5) -- (1.5,-3.5);    
    \draw [line width=2mm, \col] (-1.5,-3.5) -- (-3.5,-1.5);   
    \draw [line width=1mm] (-3.5,1.5) -- (-3.5,-1.5);
    \draw (-3.5,1.5) node (o1) {};\filldraw (o1) circle (\vsize);
    \draw (-1.5,3.5) node (o2) {};\filldraw (o2) circle (\vsize);
    \draw (1.5,3.5) node (o3) {};\filldraw (o3) circle (\vsize);
    \draw (3.5,1.5) node (o4) {};\filldraw (o4) circle (\vsize);
    \draw (3.5,-1.5) node (o5) {};\filldraw (o5) circle (\vsize);
    \draw (1.5,-3.5) node (o6) {};\filldraw (o6) circle (\vsize);
    \draw (-1.5,-3.5) node (o7) {};\filldraw (o7) circle (\vsize);
    \draw (-3.5,-1.5) node (o8) {};\filldraw (o8) circle (\vsize);
    \draw[black,thick]   (-3.5,1.5) -- (-1.5,3.5)  -- (1.5,3.5)  -- (3.5,1.5) -- (3.5,-1.5) -- (1.5,-3.5)  -- (-1.5,-3.5)  -- (-3.5,-1.5) -- (-3.5,1.5);
    \node[scale=1] at (-3.8,0.2) {\color{black}$e_1^c$};
    \node[scale=1] at (-2.7,-2.7) {\color{black}$e_x^c$};
    \node[scale=1] at (0,-3.8) {\color{black}$e_2^c$};
    \node[scale=1] at (2.7,-2.7) {\color{black}$e_y^c$};
    \node[scale=1] at (3.8,0) {\color{black}$e_3^c$};
    \node[scale=1] at (2.7,2.7) {\color{black}$e_z^c$};
    \node[scale=1] at (0,3.8) {\color{black}$e_4^c=out_c$};
    \node[scale=1] at (-2.7,2.7) {\color{black}$e_5^c$};
    \node[scale=1] at (-3.95,-0.3) {\color{black}$in_c$};
    \draw (-5,1.5) node (in_c_11) {};\filldraw (in_c_11) circle (\vsize);
    \draw (-5,-1.5) node (in_c_12) {};\filldraw (in_c_12) circle (\vsize);
    \draw (-1.5,5) node (out_c_21) {};\filldraw (out_c_21) circle (\vsize);
    \draw (1.5,5) node (out_c_22) {};\filldraw (out_c_22) circle (\vsize);
    \draw[thick] (-3.5,1.5) ..controls(-1.8,1.5) and (-1.5,1.8) .. (-1.5, 3.5);
    \draw[thick] (-3.5,-1.5) ..controls(-0.5,-1.5) and (1.5,0.5)  .. (1.5,3.5);
    \draw (-3.22,1.5) node (n1) {}; \filldraw (n1) circle (\vsize);
    \draw (-2.95,1.51) node (n2) {};\filldraw (n2) circle (\vsize);
    \draw (-2.68,1.54) node (n3) {};\filldraw (n3) circle (\vsize);
    \draw (-2.4,1.59) node (n4) {};\filldraw (n4) circle (\vsize);
    \draw (-2.16,1.675) node (n5) {};\filldraw (n5) circle (\vsize);
    \draw (-1.94,1.79) node (n6) {};\filldraw (n6) circle (\vsize);
    \draw (-1.74,2) node (n7) {}; \filldraw (n7) circle (\vsize);
    \draw (-1.64,2.25) node (n8) {};\filldraw (n8) circle (\vsize);
    \draw (-1.57,2.5) node (n9) {};\filldraw (n9) circle (\vsize);
    \draw (-1.53,2.8) node (n10) {};\filldraw (n10) circle (\vsize);
    \draw (-1.52,3.05) node (n11) {};\filldraw (n11) circle (\vsize);
    \draw (-1.5,3.3) node (n12) {};\filldraw (n12) circle (\vsize);
    \draw (-2.85,-1.45) node (n13) {}; \filldraw (n13) circle (\vsize);
    \draw (-2.35,-1.385) node (n14) {};\filldraw (n14) circle (\vsize);
    \draw (-1.8,-1.27) node (n15) {};\filldraw (n15) circle (\vsize);
    \draw (-1.2,-1.02) node (n16) {};\filldraw (n16) circle (\vsize);
    \draw (-0.6,-0.7) node (n17) {};\filldraw (n17) circle (\vsize);
    \draw (-0.1,-0.3) node (n18) {};\filldraw (n18) circle (\vsize);
    \draw (0.3,0.1) node (n19) {}; \filldraw (n19) circle (\vsize);
    \draw (0.7,0.6) node (n20) {};\filldraw (n20) circle (\vsize);
    \draw (1.02,1.2) node (n21) {};\filldraw (n21) circle (\vsize);
    \draw (1.27,1.8) node (n22) {};\filldraw (n22) circle (\vsize);
    \draw (1.385,2.35) node (n23) {};\filldraw (n23) circle (\vsize);
    \draw (1.45,2.85) node (n24) {};\filldraw (n24) circle (\vsize);
    \draw[thick] (-3.22,1.5) -- (-2.85,-1.45); 
    \draw[thick] (-2.95,1.51) --(-2.35,-1.385); 
    \draw[thick] (-2.68,1.54)  -- (-1.8,-1.27) ; 
    \draw[thick]  (-2.4,1.59) --  (-1.2,-1.02); 
    \draw[thick]  (-2.16,1.675) -- (-0.6,-0.7) ; 
    \draw[thick]  (-1.94,1.79) --  (-0.1,-0.3); 
    \draw[thick]  (-1.74,2) --  (0.3,0.1); 
    \draw[thick]  (-1.64,2.25) -- (0.7,0.6) ; 
    \draw[thick]  (-1.57,2.5) --  (1.02,1.2); 
    \draw[thick]  (-1.53,2.8) -- (1.27,1.8) ; 
    \draw[thick]  (-1.52,3.05) -- (1.385,2.35) ; 
    \draw[thick]  (-1.5,3.3) -- (1.45,2.85) ; 
    \draw[thick]  (-3.5,1.5) .. controls (-2.8,1.8) .. (-2.16,1.675) ; 
    \draw[thick]  (-1.94,1.79) .. controls (-1.95, 2.5) .. (-1.5,3.5) ; 
    \draw[thick]  (-1.2,-1.02) -- (-1.5,-3.5); 
    \draw[thick] (-0.6,-0.7)  -- (1.5,-3.5); 
    \draw[thick]  (0.7,0.6) -- (3.5,-1.5); 
    \draw[thick]  (1.02,1.2) -- (3.5,1.5); 
    \draw[thick]  (-5,1.5) --  (-3.5,1.5); 
    \draw[thick]  (-5,-1.5) -- (-3.5,-1.5) ; 
    \draw[thick]  (-1.5,3.5) -- (-1.5,5); 
    \draw[thick]  (1.5,3.5) -- (1.5,5); 
    \end{tikzpicture}
        }
        
        \end{center}
    \caption{Clause gadget for type 1 (left) and type 2 (right)}
    \label{types_of_clause_gadgets}
\end{figure}
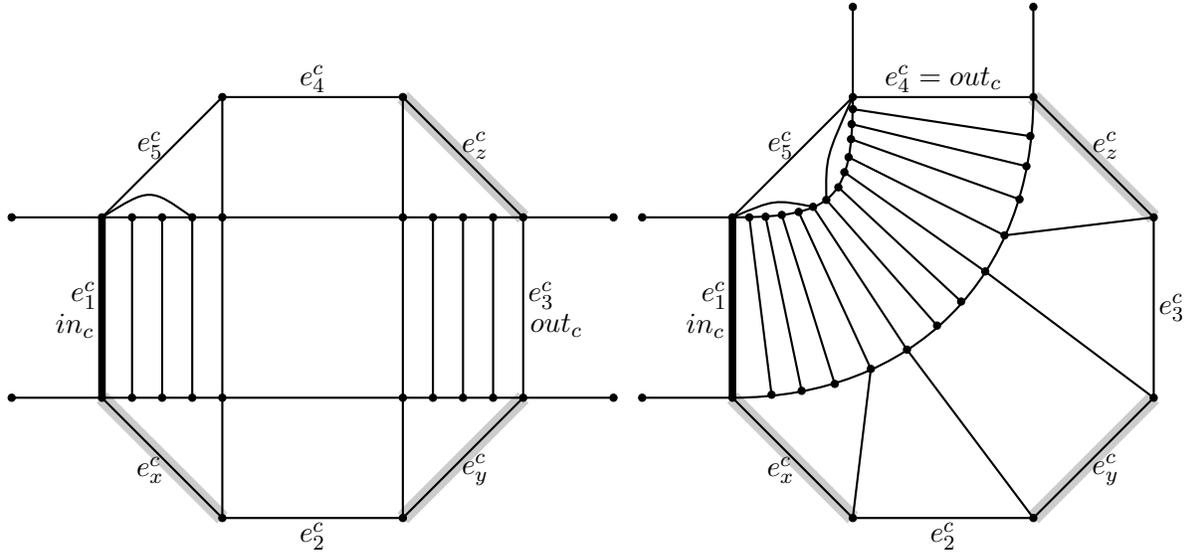

	\noindent\textbf{Clause gadget.} 
	In the planar embedding $\Pi$, each clause node $c$ is adjacent to three variable nodes, say $x,y,z$, and two more clause nodes, say $c_1$ and $c_2$. There are two types in which $c$ can be surrounded by nodes  $x,y,z,c_1,c_2$ in $\Pi$:
	only two variables appear together on one side of the clause-cycle,
	or all three variables appear consecutively on one side of the clause cycle,
	see \cref{types_of_clause_nodes}.
	For a replacement to these two types of configurations, we define two different clause gadgets (see \cref{types_of_clause_gadgets}). 
	For $c\in C$, the basic structure of the clause gadget $G_c$ includes an 8-cycle
	which contains 3 edges $e_x^c, e_y^c,$ and $e_z^c$ (corresponding to $x, y,$ and $z$), and five generic edges
	$e_1^c, e_2^c, e_3^c, e_4^c$, $e_5^c$.
	In addition, the universal ladder $\ladder_U$ passes through the gadget as shown in \cref{types_of_clause_gadgets}.
	For each of the two types of clause gadgets, we distinguish two edges of the clause gadget from where the universal ladder enters and exits by the notations $in_{c}$ and $out_{c}$ respectively.
	
	\medskip
	\noindent\textbf{Variable gadget.} Each variable node $x$ in $\Pi$ is replaced by a sufficiently large ladder-cycle $\laddercycle_x$, %
	see \cref{variable_gadget}. The bold edges in \cref{variable_gadget} correspond to the $0^{th}$-steps of $\laddercycle_x$. The orientation of all $0^{th}$-steps in $\laddercycle_x$ can either be clockwise or anticlockwise in an sc-orientation of the graph. This is because of the cyclic ladder formed in the gadget. The two types of orientations are in one-to-one correspondence with the two boolean assignments of the variable in the \tsat\ formula.
	
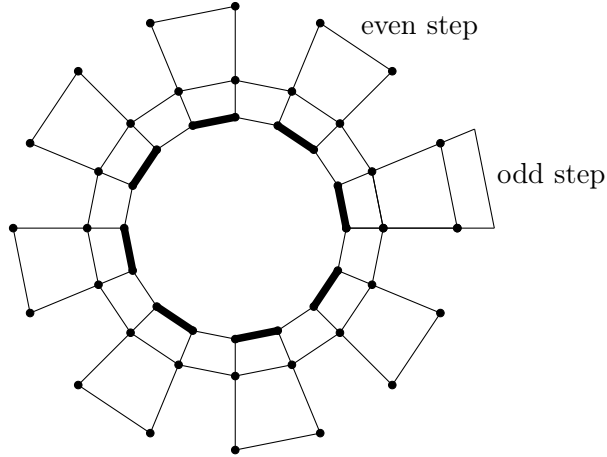
\begin{figure}[t]
    \begin{center}
    \begin{tikzpicture}[scale=0.7, node/.style={draw, circle, fill, minimum size=.1cm,inner sep=0pt}]
    
    \def\col{purple!20}
    \def\x{60}
    \def\y{80}
    \def\z{120}
    \foreach \i in {0,2,4,6,8,10,12,14}
    \draw[line width=1mm] (360/16*\i:\x pt) -- ({360/16*(\i+1)}:\x pt);
    \foreach \i in {0,...,16}
    	\draw (360/16*\i:\x pt) -- ({360/16*(\i+1)}:\x pt);
    \foreach \i in {0,...,16}
    	\draw (360/16*\i:\y pt) -- ({360/16*(\i+1)}:\y pt);
    \foreach \i in {0,...,16}{
    	\node[node] (a1) at (360/16*\i:\y pt) {};
    	\node[node] (a1) at (360/16*\i:\x pt) {};    	
    	\node[node] (a1) at (360/16*\i:\z pt) {};
    }
    \foreach \i in {0,2,4,6,8,10,12,14}
    \draw (360/16*\i:\z pt) -- ({360/16*(\i+1)}:\z pt);
    \foreach \i in {0,...,16}
    \draw (360/16*\i:\x pt) -- ({360/16*\i}:\z pt);
    
    \draw (360/16*1:\z pt) -- ({360/16*1}:140 pt);
    \draw (360/16*0:140 pt) -- ({360/16*1}:140 pt);
    \draw (360/16*0:\z pt) -- ({360/16*0}:140 pt);

    \node[scale=1] at (3.5,3.8) {$\text{even step}$};
    \node[scale=1] at (6,1) {$\text{odd step}$};

    \end{tikzpicture}
    \caption{Variable gadget (the laddercycle $\laddercycle_x$)}
    \label{variable_gadget}
    \end{center}
\end{figure}

	\medskip
	\noindent\textbf{Construction.}
	We make the following replacements while constructing our input graph $G$ of \OSC\ from the planar embedding $\Pi$ given in the input of the \cltsat\ problem:
	\begin{itemize} 
		\item Replace each clause node $c$ in $\Pi$ by a clause gadget $G_c$ in $G$. 
		\item Replace each variable node $x$ in $\Pi$ by a variable gadget $\laddercycle_x$ in $G$.  
		\item Replace each edge $\{x,c\}$ in $\Pi$ %
		by extending a ladder from the variable gadget $\laddercycle_x$ to the edge $e_x^c$ of the clause gadget $G_c$.
		If $\bar{x}$ appears in $C$, then we identify an odd step of $\laddercycle_x$ to %
		$e_x^c$. %
		If $x$ appears in $c$, then we identify an even step of $\laddercycle_x$ to the corresponding edge $e_x^c$. %
		\item For a $\{c_i,c_j\}$ edge in $\Pi$, where $c_i$ appears before $c_j$ in $\eta$, we identify the dangling edges next to $out_{c_i}$ to that of the dangling edges next to $in_{c_j}$ while maintaining the planarity.
		Moreover, let us fix the $0^{th}$-step of the universal ladder $\ladder_U$ to be the $in_{c_1}$ edge of the clause gadget $G_{c_1}$.
	\end{itemize} 
	
	Note that all clause and variable gadgets are planar, and the edge replacements by ladder also maintains the planarity of the resulting graph. 
	As per our construction of the clause gadgets, for each structure $G_c$, the $in_{c}$ and $out_{c}$ edges together are  either oriented clockwise or anti-clockwise in every sc-orientation of $G_\eta(\phi)$. Moreover, all the $in_c$ edges are oriented in one direction for all $c\in C$, and all the $out_c$ edges are oriented in the other direction. If we fix the orientation of the $0^{th}$-step of $\ladder_U$, it fixes the orientation of all edges associated with $\ladder_U$ and consequently fixes the orientation of the $in_{c}$ and $out_{c}$ edges for all $c\in C$. It also fixes the orientation of all generic edges in $G_c$ for all $c\in C$.
	We remain to prove that the \cltsat\ instance is a yes-instance if and only if the corresponding constructed instance of the \osc\ problem is a yes-instance. 
	
	\noindent\textit{Correctness of the \np-hardness reduction.} 
	As per our construction of the clause gadgets, for each structure $G_c$, the $in_{c}$ and $out_{c}$ edges together are  either oriented clockwise or anti-clockwise in every sc-orientation of $G_\eta(\phi)$. Moreover, all the $in_c$ edges are oriented in one direction for all $c\in C$, and all the $out_c$ edges are oriented in the other direction. If we fix the orientation of the $0^{th}$-step of $\ladder_U$, it fixes the orientation of all edges associated with $\ladder_U$ and consequently fixes the orientation of the $in_{c}$ and $out_{c}$ edges for all $c\in C$. It also fixes the orientation of all generic edges in $G_c$ for all $c\in C$.
	We now prove that the \cltsat\ instance is a yes-instance if and only if the corresponding constructed instance of the \osc\ problem is a yes-instance.
	
	($\Rightarrow$)  Let $\Gamma: X\rightarrow \{0,1\}$ be an assignment of variables in $X$ such that $\phi$ is satisfied. Let us fix the orientation of the $0^{th}$-step of $\ladder_U$ such that for $c\in C$, the $(in_c,out_c)$ pairs are oriented anticlockwise.
	We orient the $0^{th}$-steps of $\laddercycle_x$ anticlockwise if $\Gamma(x)=1$, and clockwise if $\Gamma(x)=0$. We claim the resulting digraph $D$ thus formed is singly connected.
	Note that, according to our construction, for an sc-orientation, the $0^{th}$-step of $\ladder_U$, and the $0^{th}$-steps of $\laddercycle_x, x\in X$, influence the orientation of all other edges in $G_{\eta}(\phi)$. Thus, the only possibilty of forming two paths between a pair of vertices is in the clause gadgets-- through paths $P_1=e^c_5$ and $P_2=e_1^c,e_x^c, e_2^c,e_y^c, e_3^c, e_z^c, e_4^c$. 
	To avoid path $P_2$, at least one edge in $\{e_x^c,e_y^c,e_z^c\}$ has to be oriented in the reverse direction to that of the generic edges in the path. This always happens because one of $\{x,y,z\}$, say $x$ satisfies $c$, and hence according to the construction $e_x^c$ gets the reverse direction.
	All other undirected cycles in $G$ contain at least 2 ladder vertices each with 2 ladder edges incident on it. Thus they are either sinks or sources in the cycle.   Hence, due to Lemma \ref{lemma:no:directed:cycles}, we know that there are no two directed paths forming in this cycle. 
	This results into an sc-orientation $D$ of $G$.

	($\Leftarrow$) If $G$ is a yes-instance of \osc, then for each clause gadget, there is at least one variable edge, say $e_x^c$, that is oriented in the  opposite direction to that of the edges $e_1^c,e_2^c,e_3^c,e_4^c$; this avoids formation of two disjoint paths between a pair of vertices in the clause gadget. The orientation of all variable edges, $e_x^c,e_y^c,e_z^c$, depend on their respective variable-gadgets. In the sc-orientation of $G$, consider the orientation of the universal ladder, that is the orientation of ($in_{c_1}, out_{c_1}$) pair, %
	and the orienations of the $0^{th}$-steps of the variable-cycles. If the orientations of $\ladder_U$ and $\laddercycle_x$ are the same, then assign value $\Gamma'(x)=1$, otherwise assign $\Gamma'(x)=0$. We claim, the resulting assignment $\Gamma'$ is a satisfying assignment of $\phi$.
	Note that each variable is either set to $1$ or $0$ in $\Gamma'$. Moreover, for each clause $c\in C$, there is at least one literal $x\in c$ which avoids the formation of paths $P_1$ and $P_2$. The corresponding literal which prevents this multi-paths is set to $1$ in $\Gamma'$, and hence each clause of $\phi$ is satisfied by assignment $\Gamma'$.
\end{prooff}

\section{Girth and Coloring Number}
\label{section:girth}

Notably the \np-hardness reduction so far heavily relied on a cycle of length $4$, as part of a `ladder'.
Hence the constructed graphs have girth $4$. We recall the \define{girth} of a graph is the length of its shortest cycle.
Can we extend this result to girth $5$ graphs or for graphs of even higher girth?

\subsection{Upper Bounds}

The \np-hardness proof cannot be extended to planar graphs of girth $5$,
since such graphs are always sc-orientable.

\newcommand{\lemmaPlanarGraphsGirthFive}{Planar graphs with girth at least 5 are sc-orientable.}
\begin{lemma}%
	\label{lemma:planar:graphs:of:girth:5}
	\lemmaPlanarGraphsGirthFive
\end{lemma}
\begin{prooff}
	We show an even stronger statement: Every near-bipartite graph $G$ of girth at least 5 is sc-orientable.
	A graph is \define{near-bipartite} if it has an independent feedback vertex set;
	that is $V(G)$ can be partitioned into $I,F$ such that $I$ is independent and $F$ is a forest.
	Then the claim about planar graphs follows due to the fact that planar graphs of girth at least 5 are near-bipartite \cite{BorGle01}.
	
	Let $I,F$ be a partition of $V(G)$ such that $G[I]$ is independent and $G[F]$ is a forest.
	Orient all edges that are incident to $I$ away from $I$; hence all the vertices in $I$ are sources.
	Fix a coloring of $G[F]$ with colors $\{1,2\}$.
	Orient every edge in $G[F]$ towards the vertex with color 2;
	hence all vertices in color class $2$ are sinks.
	Let the directed graph thus formed be~$D$.
	
	We claim that $D$ is singly connected. Assuming otherwise, let the vertices $s,t \in V(D)$ be connected by two disjoint paths $P_1,P_2$. Since $P_1,P_2$ form a cycle in the underlying graph and $G$ has girth 5, at least one of $P_1, P_2$, say $P_1$, has length at least 3.
	Let $s, v_1, v_2, v_3$ be the first 4 vertices of $P_1$. Note that $v_1, v_2\notin I$,  because both have an incoming edge. 
	Thus exactly one of $v_1$ or $v_2$ belongs to color class 2 in $G[F]$. Hence one of $v_1, v_2$ is a sink, contradicting to the existence of path $P_1$. Thus $D$ is singly connected, and $G$ is sc-orientable.
\end{prooff}

Also, for general graphs, if the girth is at least twice the coloring number, the graph is sc-orientable.

\begin{theorem}
	\label{lemma:girth:coloring}
	Let $\girth(G) \geq 2\chi(G)-1$.
	Then $G$ is sc-orientable.
\end{theorem}
\begin{prooff}
	Consider a proper vertex-coloring $c: V(G) \to \{1,\dots,|\chi(G)|\}$.
	Let $\sigma$ be an orientation which orients each edge $\{u,v\}$ to the vertex which corresponds to the higher color class, that is, to $\argmax_{w \in \{u,v\}} c(w) $.
	We claim that $G_\sigma$ is singly connected.
	Assuming the contrary, there are $s,t \in V(G)$ with disjoint $s,t$-paths $P_1,P_2$.
	Then $P_1,P_2$ form a cycle of length at least $\girth(G) \geq 2|\chi(G)|-1$.
	At least one of the paths, say $P_1$, consists of $>\chi(G)$ vertices.
	Then there are vertices $u,v$ on path $P_1$ with the same color.
	Without loss of generality say $u$ appears before $v$ in $P_1$.
	To form an $s,t$-path the successor of $u$ on path $P_1$ must have a higher color than $u$ and so on.
	Thus a later vertex $v$ with the same color contradicts that $P_1$ is an $s,t$-path.
\end{prooff}

As a remark, by Grötzsch's Theorem every planar graph of girth at least $4$ is $3$-colorable~\cite{groetzsch1959}. 
Hence, \Cref{lemma:planar:graphs:of:girth:5} also follows from \Cref{lemma:girth:coloring}.

\subsection{The Dichotomy}
\label{subsection:dichotomy}

We do not know whether there even exists a graph of girth at least $5$ that is not sc-orientable.
Intriguingly, however, we do not need to answer this question to follow a dichotomy between the \p-time and \np-completeness.
That is, for every $g,c\geq 3$, $\textsc{SC-Orienation}$ when restricted to graphs of girth $g$ and coloring number $c$ is either \np-complete or trivially solvable. Let the class of graphs of girth $g$ and chromatic number $c$ be $\GG_{g,c}$.  
Our proof does not pinpoint the exact boundary between \p-time and \np-completeness.
We only know, for a fixed $c$, there is some upper bound $g'$ for the values of $g$ such that the hardness is only unclear on $\GG_{g,c}$, where $g<g'$. See \Cref{lemma:girth:coloring}. Moreover, we showed in \Cref{section:planar}, the \np-hardness on $\GG_{g,c}$ for $g=3,4$.
The key ingredient for our dichotomy result is to compile an \np-hardness proof that only relies upon a single no-instance of a certain girth $g$ and coloring number $c$.

The main building block of our construction is a \emph{coupling gadget}.
It is a graph that couples the orientation of two special edges
but that does \emph{not} enforce a coloring depending on the orientation.
An sc-orientable graph $H$ with special edges $\{a_1,b_1\}, \{a_2,b_2\} \in E(H)$
is a \define{coupling gadget} on $(a_1,b_1), (a_2,b_2)$
if
\begin{enumerate}[{(1)}]
	\item \sloppy
	$(a_1,b_1), (a_2,b_2)$ are \define{coupled},
	that is, for every sc-orientation $\sigma$,
	either $\sigma(\{a_x,b_x\}) = b_x$ for $x\in\{1,2\}$,
	or $\sigma(\{a_x,b_x\}) = a_x$  for $x\in\{1,2\}$, \label{property1}
	\item there is an sc-orientation without any directed path from $\{a_1,b_1\}$ to $\{a_2,b_2\}$,
	and\label{property2}
	\item there are two $\chi(H)$-colorings $f,f'$ of $H$
	such that
	$f(a_1)=f(a_2)$, $f(b_1)=f(b_2)$
	and 
	$f'(a_1)=f'(b_2)$, $f'(b_1)=f'(a_2)$.\label{property3}
\end{enumerate}

For example, $G_{2,3}$ ladder (which we use for \Cref{lemma:np:planar}) satisfies only properties (1) and~(2).
In turn, $G_{2,4} \cup K_3$ is a coupling gadget.

\begin{lemma}
	\label{lemma:sc:g:c}
	For every $g,c \geq 3$, there is an sc-orientable graph $G \in \GG_{g,c}$.
\end{lemma}
\begin{prooff}
	We know that there is at least one graph $G_{g',c}$ of girth $g'$ and chromatic number~$c$, see~\cite{erdos1959graph}.
	When $g=g' \geq 2c-1$, \Cref{lemma:girth:coloring} yields that $G_{g,c}$ is sc-orientable.
	For $g=g' < 2c-1$, the graph, $G_{2c-1,c} \cup C_g,$ is sc-orientable; it has girth $g$ and chromatic number $c$.
\end{prooff}

\newcommand{\lemmaNoIffCoupling}{
	For $g\geq3$, $c\geq4$, %
	there is a non-sc-orientable graph $G \in \GG_{g,c}$,
	if and only if
	there is a coupling gadget $H \in \GG_{g,c}$.
}

\begin{lemma}%
	\label{lemma:no:iff:coupling}
	\lemmaNoIffCoupling
\end{lemma}
\begin{prooff}
	\forward
	Consider a non-sc-orientable graph $G \in G_{g,c}$.
	We want to assume that $G$ contains at least one edge $\{u,v\}$
	such that subdividing $\{u,v\}$ yields an sc-orientable graph.
	If that is not the case, pick an edge $\{u,v\}$	and consider the graph $G'$ resulting from graph $G$ where $\{u,v\}$ is subdivided into a path $uwv$.
	To assure termination, initially assume all edges unmarked,
	and consider subdividing an unmarked edge $\{u,v\}$ before considering any marked edges.
	Whenever an edge $\{u,v\}$ is subdivided,
	mark the new edges $\{u,w\}$ and $\{w,v\}$.
	The procedure terminates at the latest when all edges are marked
	since then the resulting graph is bipartite
	and hence trivially sc-orientable.
	Each above step does not decrease the girth and does not increase the coloring number.
	
	Now, let $H'$ consists of a copy of $G$ where edge $\{u,v\}$ is subdivided resulting in a path $uwv$ such that  $H'$ has an sc-orientation $\sigma$.
	Note that $girth(H')\geq girth(G)$.
	Crucially, $\sigma$ couples edges $(u,w)$ and $(w,v)$,
	meaning it orientates $\{u,w\}, \{w,v\}$ either both away from $w$ or both towards $w$.
	To see this, assume an sc-orientation $\sigma'$ that, up to symmetry,
	has $\sigma'(\{u,w\}) = w$ and $\sigma'(\{w,v\}) = v$.
	Then $\sigma'$ extended with $\sigma'(\{u,v\}) = v$ forms an sc-orientation of $G$,
	despite it being non-sc-orientable.
	Thus indeed edges $(u,w),(v,w)$ are coupled.
	
	To construct $H$, we introduce four copies $H_1,H_2,H_3,H_4$ of $H'$ naming their vertices with subscript $1,2,3,4$ respectively.
	We identify the vertices $w_1,u_2$ and vertices $v_1,w_2,u_3$ and vertices $v_2,w_3,u_4$ and vertices $v_3,w_4$; see the gadget in \cref{figure:small:graphs}.
	(Eventually our graph $H$ results from adding another sc-orientable graph from $\GG_{g,c}$, which we can ignore for now.)
	By the transitivity of the coupled edges in $H_1,H_2,H_3,H_4$,
	we have $(u_1,w_1),(w_4,v_4)$ coupled, showing Property \ref{property1}.
	In particular, combining sc-orientations of $H_1,H_2,H_3,H_4$ that agree on the orientation of the edges on the path
	$P = u_1,u_2,u_3,u_4,w_4,v_4$ (which have to be alternately orientated), yields existence of an sc-orientation.
	To show Property \ref{property2}, consider any sc-orientation $\sigma$.
	Up to symmetry, $\sigma$ has the edges on $P$ oriented towards $w_1$ and $w_3$.
	Then $G_\sigma$ contains no directed path from $u_1$ to $v_1$,
	since this would imply two distinct $u_1,w_1$-paths.
	Similarly, there is no directed path between $u_2,v_2$,
	and we conclude that $H$ satisfies Property \ref{property2}.
	
	Finally, we state colorings $f$ and $f'$ with colors $\{1,\dots,\chi(H')\}$ that show Property \ref{property3}.
	It suffices to define the coloring of the vertices on the path $P$,
	as long as each three consecutive vertices have pairwise distinct colors.
	Then each remaining un-colored component $C$ is attached to two vertices $x,y \in V(P)$ in distance two on $P$,
	where $G[C \cup \{x,y\}]$ is a copy of $H'$
	and hence can be $\chi(G)$-colored when the color of $x,y$ is already fixed.
	To obtain a coloring $f$, let $P$ be colored $1,2,3,4,1,2$, and in order 
	to obtain a coloring $f'$, let $P$ be colored $1,2,3,4,2,1$.
	Here we use that $c\geq 4$.
	
	So far, the coloring number did not increase as the coloring $f$ and $f'$ demonstrate.
	Also the girth does not decrease
	since every chordless cycle is %
	contained in a copy of $H'$.
	In order to ensure $H \in \GG_{g,c}$, for the final construction step, let $H$ be the disjoint union of two graphs,
	$H$ and $G_{g,c}$, where $G_{g,c} \in \GG_{g,c}$.
	Since $G_{g,c}$ forms its own component, it does not interfere with the properties (1)-(3) observed so far.
	Thus $H$ is a coupling gadget on $(u_1,w_1)$ and $(w_4,v_4)$.
	
	\backward
	Let $H \in \GG_{g,c}$ be a coupling gadget on $(a_1,b_1)$ and $(a_2,b_2)$.
	Let $g' \in \{g,g+1\}$ be odd.
	We construct $G$ from $g'$ copies of $H$, denoted as $H^1,\dots,H^{g'}$.
	Essentially we glue the coupling gadgets together at their coupled edges in a cycle
	while introducing a final `twist' of the coupling edges.
	The graph $G$ result from identifying $a^i_2, a^{i+1}_1$ and identifying $b^{i}_2,b^{i+1}_1$ for $i \in \{1,\dots, g'-1\}$,
	and identifying $a^{g'}_2, b^{1}_1$ and identifying $b^{g'}_2,a^1_1$.
	We further add an sc-orientable graph $G_{g,c} \in \GG_{g,c}$ to our construction. Such a graph always exists because of \Cref{lemma:sc:g:c}.
	Since $g'$ is odd, not all $(a_1^i,b_1^i)$ and $(a_2^i,b_2^i)$ for $i \in \{1,\dots g'\}$ can be coupled to obtain an sc-orientation of $G$.
	Thus $G$ is not sc-orientable.
	
	Graph $G$ has the same girth as $H$,
	since any chordless cycle either has length $\geq g'$ or is contained in a copy of $H$.
	Also we can color $G $ with $c$ colors using the $c$-colorings $f,f'$ of Property \ref{property3}.
	We use $f$ to color $H^i$ with $f(a_1^i)=f(a_2^i)$ and $f(b_1^i)=f(b_2^i)$, for $i \in \{1, \dots, g'-1\}$.
	Finally, we use the $c$-coloring $f'$ to color $H^{g'}$
	which then agrees on the coloring of the vertices shared with $H^1$ and shared with $H^{g'-1}$.
\end{prooff}

\newcommand{\lemmaDicho}{
	For $g,c\geq 3$,
	\textsc{SC-Orientation} restricted to the graphs in $\GG_{g,c}$ is either \np-complete or trivially solvable.
}
\begin{theorem}%
	\label{lemma:dicho}
	\lemmaDicho
\end{theorem}

\begin{prooff}
	Consider fixed $g,c\geq 3$. If all graphs in $\GG_{g,c}$ are sc-orientable, then an algorithm may always answer `yes' to decide whether they are sc-orientable. 
	If not, there exists at least one no-instance of \osc\ in $\GG_{g,c}$. In this case, 
	we claim, it is \np-hard to decide whether an arbitrary graph in $\GG_{g,c}$ is sc-orientable. 
	\Cref{lemma:no:iff:coupling} derives a coupling gadget $H \in \GG_{g,c}$ when the problem is not trivial on $\GG_{g,c}$.
	
	In order to prove the \np-hardness claim, we draw ideas from the \np-hardness construction given in \Cref{lemma:np:planar}. Recall, in \Cref{lemma:np:planar}, each clause gadget consists of a cycle which has generic edges and variable edges. 
	Moreover, the orientation of generic-edges and variable-edges are fixed once we fix the orientation of the $0^{th}$-steps of the universal ladder and the variable-gadgets respectively. Let $\sigma$ be such a fixed orientation.  Note that, $\sigma$ is not an sc-orientation if and only if there is a clause gadget whose all variable-edges are oriented in the same direction as that of edge $e_1^c$, and thus forming two disjoint paths in the cycle. We exploit this property of the clause gadget. \\
	
	\noindent\textbf{Adaptation to the \np-hardness construction in \Cref{lemma:np:planar}:}
	\begin{itemize}
		\item
		We construct an even length clause-gadget of size at least $g$.
		Let $g' \in \{g,g+1\}$ be odd . %
		Update the construction of clause-gadget in \Cref{lemma:np:planar}, and replace the edge $e_2^c$ by a set of edges $\{e_{2_1}^c, e_{2_2}^c, \dots, e_{2_{g'}}^c\}$.
		The length of the clause-cycle in the clause-gadget is at least $g$.
		Moreover, there exists a 2-coloring of the vertices of each clause-cycle. 
		\item The universal ladder and the variable-gadgets (in this case, the variable-ladders) each consists of plenty of copies of the coupling gadgets $H\in \GG_{g,c}$ stacked on each other such that each pair of consecutive coupling gadgets share a unique coupled edge.
		Let us refer these unique edges as the steps of the ladder.
		\item According to Property \ref{property2}, the coupling-edges $(a_1, b_1)$ and $(a_2,b_2)$ of a coupling-gadget are disjoint. Thus one can introduce a coupling gadget between a pair of disjoint edges to make them coupled.
		Note that for an arbitrary edge $\{a,b\}$, we can couple $(a',b')$ with $(a,b)$, or $(b,a)$.
		The latter is simply done by introducing a twist in the graphical representation of the structure. Moreover, to avoid short cycles of length less than $g$, we can do the coupling by introducing a stack of $g$ many coupling gadgets instead of one coupling gadget between a pair of disjoint edges.
		Let us say $(a', b')$ is coupled with $(a,b)$, and $(a',b')$ is reverse-coupled with $(b,a)$.
		Using this idea, we connect all the generic edges of the clause-gadgets to the universal ladder.
		We couple $e_{1}^c, e_{3}^c, e_{4}^c$ and $e_{2_1}^c, e_{2_2}^c, \dots, e_{2_{g'}}^c$, and reverse-couple the edge $e_{5}^c$ with the steps of the universal ladder.
		\item Similarly, we couple variable-edges with the steps of their respective variable ladders.
		We couple $e_x^c$ with the variable ladder if $x\in c$, else we reverse-couple $e_x^c$ with the variable ladder if $\bar{x}\in c$.
	\end{itemize}
	
	\begin{sloppypar}
		\noindent\textbf{Correctness of the construction:}
		The constructed instance of the \osc\ problem indeed has girth $g$. This is because, since the coupling gadget $H$ belongs to class $\GG_{g,c}$, the universal ladder and the variable gadget have girth $g$, and the clause-cycle has length at least $g$. In addition we maintain distance at least $g$ between an edge of the clause-gadget and an edge of a (universal/variable) ladder. 
		We claim that the constructed instance has chromatic number $c$.
		Each clause-gadget can be colored with colors $\{1,2\}$.
		We color each step of the ladders also with colors $\{1,2\}$.
		It remains to extend this coloring to the internal vertices of each coupling-gadget.
		This is possible because of Property \ref{property3} of the coupling gadgets: There exists two $\chi(H)$-colorings $f,f'$ of the coupling gadget $H$ such that $f(a_1)=f(a_2)$, $f(b_1)=f(b_2)$ and  $f'(a_1)=f'(b_2)$, $f'(b_1)=f'(a_2)$.
		Hence each coupling-gadget has a coloring using $f$ or $f'$ that respects the prescribed coloring.
		Also, note that the coupling gadget $H$ itself cannot be colored with less than $c$ colors, hence the chromatic number of the entire graph is $c$. 
	\end{sloppypar}
	
	The proof of \Cref{lemma:np:planar} has been adapted in this settings in a way that the rest of the arguments about the correctness of the construction follow as it is. Thus we conclude our proof here. 
\end{prooff}

\commentt{
	\begin{prooff}
		\todo{original}
		Consider fixed $g,c\geq 3$.
		If all graphs $G \in \GG_{g,c}$ are sc-orientable,
		then an algorithm may simply answer always `yes'.
		Otherwise, \Cref{lemma:no:iff:coupling} provides a coupling gadget $H \in \GG_{g,c}$.
		We can adapt the \np-hardness construction of \Cref{lemma:np:planar}.
		First, we avoid ladders $L$ of an odd number of $C_4$
		(attached to a clause gadget where the corresponding variable is used negated)
		by essentially placing another $C_4$ between the end of the latter and the corresponding edge in th clause gadget:
		Instead of identifying $(u,v)$ of the ladder with $(u',v')$ of a clause gadget,
		add edges $\{u,v'\}$ and $\{v,u'\}$.
		(\Cref{lemma:np:planar} avoids such a construction to yield a planar graph.)
		Hence we can assume that each ladder is made up of the coupling gadget of girth $4$, namely a domino (with a triangle as an extra component).
		We replace each such coupling gadget by a coupling gadget $H \in \GG_{g,c}$.
		We may insert $g$ coupling gadgets between each original domino to have any cycle using the ladders length at least $g$.
		Similarly, we can enlarge the ladders within the clause gadget to avoid short cycles.
	\end{prooff}
}

\section{Perfect Graphs}
\label{section:perfect}

This section considers perfect graphs.
First we show that the \np-hardness proof from \Cref{lemma:np:planar} can be adapted to perfect graphs.
Then later we study the distance-hereditary graphs as a special case.

We modify our \np-hardness reduction from \Cref{lemma:np:planar} as follows:
In the reduced instance, for each edge $e=\{u,v\}$, introduce two new vertices $w_e$ and $w_e'$.
Delete edge $e$, and instead add edges such that there is a path $u w_e v$, and a triangle $u w_e w_e'$.

Then according to \Cref{lemma:triangle:shrink}, the modified instance is equivalent to the original instance in the sense that the modified graph is sc-orientable if and only if the originally constructed graph is sc-orientable.
Further, the modified graph is 3-colorable and has a maximum clique of size 3;
hence it is perfect.
From this we conclude the following result.

\begin{theorem}
	\textsc{SC-Orientation} is \np-hard even for perfect graphs.
\end{theorem}

\subsection{Distance-Hereditary Graphs}

Now we derive two classifications of the distance hereditary graphs that are sc-orientable.
One way is to simply restrict the recursive definition to (locally) avoid a diamond subgraph.
The second is a concise classification by forbidden subgraphs.

A graph is \define{distance hereditary} if it can be constructed starting from a single isolated vertex with the following operations~\cite{BandeltM86}:
a \emph{pendant} vertex to $u \in V(G)$, that is, add a vertex $v$ and edge $\{u,v\}$;
\emph{true twin} on a vertex $u \in V(G)$, that is, add a vertex $v$ with neighborhood $N(v)=N(u)$ and edge $\{u,v\}$;
\emph{false twin} on a vertex $u \in V(G)$, that is, add a vertex $v$ with neighborhood $N(v)=N(u)$.

We denote a graph as \define{strongly distance hereditary} if it can be constructed from a single isolated vertex with the following restricted operations:
\begin{itemize}
	\item
	a pendant vertex to $u \in V(G)$;
	\item
	true twin on a vertex $u \in V(G)$, restricted to $u$ where $|N(u)|=1$; and
	\item
	false twin on a vertex $u \in V(G)$, restricted to $u$ where $N(u)^2 \cap E(G) = \emptyset$.
\end{itemize}

Observe that the forbidden operations would immediately imply a diamond as a subgraph.
Thus every distance hereditary graph that has an sc-orientation must be strongly distance hereditary.
Now, we show also the converse is true.

\begin{theorem}%
	\label{lemma:pseudo:distance:hereditary}
	\lemmaDHarePDH
\end{theorem}
\begin{prooff}
	As discussed, it remains to show
	that every strongly distance hereditary graph indeed allows an sc-orientation.
	We show the following property for strongly distance hereditary graph~$G$:
	\begin{itemize}
		\item Every 2-vertex connected component is bipartite or a triangle.
	\end{itemize}
	This property implies that $G$ is sc-orientable,
	since $G$ decomposes into bipartite graphs and triangles, which are positive base cases of our decomposition.
	Checking this property is possible polynomial time.
	To specifically check for strongly distance hereditaryness, an additional test of whether $G$ is distance hereditary suffices; this is possible in linear time \cite{DBLP:journals/tcs/DamiandHP01,DBLP:journals/dam/HammerM90}.

	It remains to show the above property.
	An isolated vertex is a bipartite graph.
	We proceed inductively:
	(Let $G'$ below be the graph resulting from an operation)
	\begin{itemize}
		\item
		A pendant $v$ to vertex $u \in V(G)$:
		The 2-vertex connected components of $G'$ are those of $G$ and $\{v\}$.
		Component $\{v\}$ is bipartite, and components of $G$ remain bipartite or a triangle.
		\item
		True twins, for a vertex $u \in V(G)$ if $|N(u)|=1$:
		Let $N(u)=\{w\}$.
		The 2-vertex connected components of $G'$ are $\{u,v,w\}$ and those of $G$.
		Component $\{u,v,w\}$ is a triangle, and components of $G$ remain bipartite or a triangle.
		\item
		False twins, for a vertex $u \in V(G)$, where $N(u)^2 \cap E(G) = \emptyset$:
		Let $\mathcal{C}$ be the 2-vertex connected components of $G$.
		Let $u$ be a part of a 2-vertex connected component $C \in \mathcal{C}$.
		Note that $C$ is not a triangle since $N(u)^2 \cap E(G) = \emptyset$.
		Thus $C$ is bipartite.
		Then the 2-vertex connected components of $G'$ are $C \cup \{v\}$ and components $\mathcal{C} \setminus C$.
		The component $C \cup \{v\}$ is still bipartite; the vertex $u$ and $v$ can be put together in the same part of the bipartition. 
		This satisfies the property,
	\end{itemize}
	This completes the proof.
\end{prooff}

Note that distance hereditary graphs are exactly the (house, hole, domino, gem)-free graphs \cite{BandeltM86},
where a hole is any cycle of length $\geq 5$.
When we replace `gem' with `diamond' we exactly end up with the strongly distance hereditary graphs.

\newcommand{\lemmaPDHcharact}{
	The strongly distance hereditary graphs are exactly the graphs with house, hole, domino and diamond as forbidden subgraph.
}

\begin{lemma}%
	\label{lemma:pdh:charact}
	\lemmaPDHcharact
\end{lemma}

\begin{prooff}
	We show both inclusions:
	Let HHDD-free graphs be those graph with house, hole, domino and diamond as a forbidden subgraph.
	
	($\subseteq$)
	Let $G$ be a strongly distance hereditary graph,
	that is distance hereditary and sc-orientable; see \Cref{lemma:pseudo:distance:hereditary}.
	Then $G$ has house, hole and domino as well as %
	diamond as a forbidden subgraphs;
	hence is HHDD-free.
	
	($\supseteq$)
	Let $G$ be HHDD-free.
	Then $G$ is also gem-free and hence is distance hereditary.
	We show that every 2-vertex connected component of $G$ is bipartite or a triangle
	which, as per \Cref{lemma:pseudo:distance:hereditary}, implies that $G$ is sc-orientable.
	Consider a 2-vertex connected component $C$ of $G$ which is neither bipartite nor a triangle. Note that $C$ contains at least 3 vertices. Moreover, since $C$ is hole-free and not bipartite, it contains a triangle $uvw$.
	As $C$ is not merely a triangle, there is at least one other vertex $u^\star \in C \setminus \{u,v,w\}$. Since $C$ is a 2-vertex connected component, $u^\star$ must have paths to two vertices of $u,v,w$, say $P_1$ is a path from $u^\star$ to $u$ and $P_2$ is a path from $u^\star$ to $v$. Depending on the size of $P_1$ and $P_2$, we either get a hole, a diamond, or a house. Thus forming a contradiction to $G$ being HHDD-free.
	Therefore every 2-vertex connected component of the graph is bipartite or a triangle,
	which finishes the proof.
\end{prooff}

\section{Conclusion}
\label{section:conclusion}

Among others results,
we have shown that for every $g,c \geq 3$, the restriction to graphs of girth $g$ and chromatic number~$c$
is either \np-complete or in \p.
While we know \np-completeness for low girth and chromatic number
and we know that for relatively large girth compared to $\chi(G)$ the problem is trivial,
it yet remains to pinpoint the exact boundary of the \np-hard and \p cases.
As shown, to extend the \np-hardness result, one needs to simply find a no-instance of the girth and chromatic number of interest (or alternatively merely a coupling gadget),
for example graphs of girth $5$ and chromatic number $4$.
On the other hand, one might improve \Cref{lemma:girth:coloring} to also hold for smaller values of girth~$g$.

\bibliographystyle{plain}
\bibliography{mybibliography}

\newpage

\appendix

\end{document}